\documentclass[11pt]{article}
\usepackage{amsmath,amsthm,amsfonts,verbatim}

 \DeclareMathOperator{\lca}{lca}
 \DeclareMathOperator{\diam}{diam}

\newcommand{\e}{\varepsilon}

\theoremstyle{plain}
 \newtheorem{theorem}{Theorem}
 \newtheorem{lemma}{Lemma}

 \newtheorem{proposition}[lemma]{Proposition}

 \theoremstyle{definition}
 \newtheorem{definition}{Definition}

\begin{document}

\title{On Metric Ramsey-type Dichotomies}

\date{}

\author{Yair Bartal\thanks{Supported in part by a grant from the
Israeli National Science Foundation} \and Nathan
Linial\thanks{Supported in part by a grant from the Israeli
National Science Foundation} \and Manor Mendel\thanks{Supported in
part by the Landau Center.} \and Assaf Naor }

\maketitle

\begin{abstract}
The classical Ramsey theorem, states that every graph
contains either a large clique or a large independent set. Here
we investigate similar dichotomic phenomena in the context of
finite metric spaces. Namely, we prove statements of the
form "Every finite metric space contains a
{\em large} subspace that is {\em nearly}
equilateral or {\em far} from being equilateral".
We consider two distinct interpretations for
being "far from equilateral".
Proximity among metric spaces is quantified through
the metric distortion $\alpha$.
We provide tight asymptotic answers for these problems. In
particular, we show that a phase transition occurs at $\alpha=2$.
\end{abstract}

\section{Introduction}

A Ramsey-type theorem states that large systems necessarily
contain large, highly structured sub-systems. Here we consider
Ramsey-type problems for finite metric spaces, and interpret
``highly structured" as being embeddable with low distortion in
some ``simple" metric spaces.

A mapping between two metric spaces $f:M
\rightarrow X$, is called an embedding of $M$ in $X$.
The \emph{distortion} of the embedding is defined as
\[
\mathrm{dist}(f)=\sup_{\substack{x,y\in M\\x\neq
y}}\frac{d_X(f(x),f(y))}{d_M(x,y)}\cdot \sup_{\substack{x,y\in M\\x\neq
y}}\frac{d_M(x,y)}{d_X(f(x),f(y))}.
\]
The least distortion attainable by any embedding of $M$ in $X$ is
denoted by $c_X(M)$. When $c_X(M)\leq \alpha$ we say that $M$
$\alpha$-embeds in $X$. When $M$ $\alpha$-embeds in $X$ via a
bijection, we say that $M$ and $X$ are $\alpha$-equivalent.

This paper deals with the following notion.

\begin{definition}[Metric Ramsey function] For a given class of metric
spaces $\mathcal{X}$ we denote by $R_{\mathcal{X}}(\alpha,n)$ the
largest integer $m$ such that any $n$-point metric space has a
subspace of size $m$ that $\alpha$-embeds into some $X\in
\mathcal{X}$. When $\mathcal{X}=\{\ell_p\}$ we use the notation
$R_{p}$.
\end{definition}

In \cite{bfm}, Bourgain, Figiel, and Milman study this function for
$\mathcal{X}=\{\ell_2\}$, as a non-linear analog of the classical Dvoretzky
theorem~\cite{dvo}. They prove

\begin{theorem}[\cite{bfm}] \label{thm:bfm86}
For any $\alpha>1$ there exists $C(\alpha)>0$ such that $R_2(\alpha,n) \geq
C(\alpha) \log n$. Furthermore, there exists $\alpha_0>1$ such that
$R_2(\alpha_0,n)=O(\log n)$.
\end{theorem}

Lower bounds which improve on Theorem \ref{thm:bfm86} for {\em
large} $\alpha$ were obtained in \cite{bbm}, and the Euclidean
metric-Ramsey problem was comprehensively studied in~\cite{blmn-phenomena}.
There, the lower bound on $R_2$ was achieved via embedding into a
special type of $\ell_2$ metrics, namely \emph{ultrametrics}.
Denote by UM the class of ultrametrics. The following phase
transition was established.

\begin{theorem}[\cite{blmn-phenomena}]
\label{thm:l2-phase} Let  $n\in \mathbb{N}$. Then:
\begin{enumerate}
\item For every $1<\alpha<2$,
$$ c(\alpha)\log n\le
R_{\text{UM}}(\alpha,n) \le R_2(\alpha,n)\le 2\log n+C(\alpha),
$$
where $c(\alpha),C(\alpha)$ may depend only on $\alpha$.

\item For every $\alpha>2$,
$$
n^{c'(\alpha)}\le R_{\text{UM}}(\alpha,n) \le R_2(\alpha,n)\le
n^{C'(\alpha)}, $$
 where $c'(\alpha), C'(\alpha)$ depend only on
$\alpha$ and satisfy $\max\left\{0, 1- \frac{c \log
\alpha}{\alpha}\right\}< c'(\alpha)<C'(\alpha)< \min\left\{1,1 -
\frac{C}{\alpha}\right\}$, with $c,C>0$ universal constants.
\end{enumerate}
\end{theorem}

In \cite{blmn2}, a similar phase transition phenomenon is proved
for embeddings in $\ell_p$, $p\in[1,2)$.

\medskip
A natural refinement of ultrametrics was suggested in \cite{bartal1}.

\begin{definition}[\cite{bartal1}]\label{def:hst}
For $k\geq 1$, a $k$-\emph{hierarchically well-separated tree} ($k$-HST) is a
metric space whose elements are the leaves of a rooted tree $T$. To each
vertex $u\in T$, a label $\Delta(u) \ge 0$ is associated such that
$\Delta(u)=0$ iff $u$ is a leaf of $T$. The labels are such that if a vertex
$u$ is a child of a vertex $v$ then $\Delta(u)\leq \Delta(v)/k$ .
 The distance between two leaves $x,y\in T$ is
defined as $\Delta(\lca(x,y))$, where $\lca(x,y)$ is the least
common ancestor of $x$ and $y$ in $T$. $T$ is called \emph{the
defining tree} of the HST.

The notion of an ultrametric is easily seen to coincide with
that of a 1-HSTs.
\end{definition}

In \cite{blmn-phenomena}, the Ramsey problem of embedding into $k$-HSTs was
also studied.
\begin{theorem}[\cite{blmn-phenomena}] \label{thm:ramsey-hst}
For any $\e\in(0,1]$, and any $k\geq 1$,
 \[ R_{k\text{\textrm{-HST}}}(2+\e,n) \geq n^{\frac{c \e}{\log(2k/\e)}} . \]
\end{theorem}

\medskip

In this note we study Ramsey problems closer in spirit to the original
Ramsey problem in combinatorics, which is of a dichotomic nature. More
specifically, such theorems state that every metric space contains a large
subspace which is close to one of two extremal types of simple metric spaces.
In this note we consider two different (but related) type of dichotomies.

We begin with some motivation. Since every $3$-point metric is isometric to a
Euclidean triangle, we can associate with it three angles. We say that two
$3$-point metrics are $\epsilon$-similar, if the corresponding angles differ
by at most $\epsilon$. Fix some $\epsilon > 0$. The collection of all
$3$-point metrics can be partitioned into a constant number of classes such
that every two triples in the same class are $\epsilon$-similar. By Ramsey's
theorem, every $n$-point metric space contains a large homogeneous subset,
namely a set of $f=f(n)$ elements, every two triples in which are
$\epsilon$-similar, where $f$ tends to $\infty$ with $n$. It is not hard,
however, to show that there are only two types of unboundedly large
homogeneous sets. Either all triples in such a class are $\epsilon$-similar
to the equilateral triangle with angles $(60^{\circ}, 60^{\circ},
60^{\circ})$ or all are $\epsilon$-similar to a triple in which the smallest
angle is at most $\epsilon^2$, say. In the latter case, in fact, more is true.

\paragraph{The Metric dichotomy.}
In the first type of dichotomy treated, which we call the {\em
metric dichotomy}, we have on one hand equilateral spaces, i.e.
metric spaces in which all pairwise distances are equal. The
``opposite" extreme are spaces in which every triple of points is
far (in the sense of metric distortion) from being an equilateral
triangle. We define $F_k(\alpha,n)$ as the largest $m$ such that
any $n$ point metric space contains an $m$-point subspace which is
either $\alpha$-equivalent to an equilateral space or
$\alpha$-equivalent to a space for which every triple of points
has distortion at least $k$ from an equilateral triangle.

The notion of spaces in which no triple is $k$-equivalent to an
equilateral triangle is quite natural. It turns out, however, that
in order to analyze the behavior of $F_k$, it is more convenient
(and essentially equivalent) to consider instead binary $k$-HSTs,
i.e. $k$-HSTs whose defining tree is binary (every vertex has at
most two children). The relevant dichotomic Ramsey function is
defined as
 \[ E_k(\alpha,n)=R_{\left\{ \substack{ \text{binary $k$-HSTs}\\
\text{or equilateral spaces}}\right\}}(\alpha,n). \]

The relation between these notions is clarified in the following
proposition.

\begin{proposition} \label{prop:approx-equiv} The following two
assertions hold:

\begin{enumerate}
\item Let $M$ be a binary $k$-HST and let $S\subset M$, have
cardinality $|S|\geq 3$. Then
$$
c_{\{\mathrm{equilateral\  spaces}\}}(S)\geq k.
$$

\item Let $M$ be a metric space in which
$c_{\{\mathrm{equilateral\ spaces}\}}(S)\geq k$ for every
$S\subset M$ with $|S|\geq 3$, where $k> 2$. Then $M$ is
$\frac{k}{k-2}$-equivalent to a binary $\frac{k}{2}$-HST.
\end{enumerate}
In particular
 \[ E_k(\alpha,n)\leq F_k(\alpha,n) \leq E_{k/2}\left(\alpha \tfrac{k}{k-2}, n\right)
 .\]
\end{proposition}


\begin{theorem}[The metric dichotomy]\label{thm:second}
\mbox{\ }
\begin{enumerate}
\item For $\alpha>2$, $k>2$:
\[
\exp \left ( c(\alpha,k) \sqrt{\log n} \right) \leq E_k(\alpha,n) \leq
F_k(\alpha,n) \le \exp \left ( C(\alpha,k) \sqrt{ \log n} \right )
\]
\item For $1<\alpha<2$, $k>2$:
\[
c(\alpha,k) \cdot \frac{\log n}{\log \log n} \leq E_k(\alpha,n) \leq
F_k(\alpha,n)  \leq C(\alpha,k)  \frac{\log n}{\log \log n}.
\]
\end{enumerate}
Here $c(\alpha,k),C(\alpha,k)>0$ depend only on $\alpha$ and $k$. The bounds
above on $E_k$ also hold for $k\in (1,2)$.
\end{theorem}


This dichotomic Ramsey problem was first studied implicitly in \cite{bkrs}.
It is possible to deduce from their work that $E_{\log n}(4,n) \geq \exp(
c\sqrt{\log n / \log \log n})$.
A closely related problem was formulated in \cite{bbm}, where some bounds on
$E_k(\alpha,n)$ were given.

%

\medskip

\paragraph{The equilateral/lacunary dichotomy.}
Another type of dichotomy that we study, was first formulated in
\cite{krr}. On the one hand, we have again the equilateral metric
spaces. At the other extreme of the dichotomy is a class of metric
spaces in which the set of pairwise distances are sparse, which we
call {\em lacunary} metric spaces. Recall that the sequence
$a_1\ge a_2\ge\ldots\ge a_n>0$ is called $k$-lacunary for some
$k\ge 1$, if $a_{i+1}\le a_i/k$ for $i=1,\ldots,n-1$. A metric $d$
on $\{1,\ldots,n\}$ is called $k$-lacunary if there exists a
$k$-lacunary sequence $a_1\ge\ldots\ge a_n>0$ such that for $1\le
i<j\le n$, $d(i,j)=a_i$. Alternatively, $k$-lacunary spaces can be
defined using HSTs.

\begin{definition} Let $k>1$. A $k$-increasing metric space is a $k$-HST $X$ such that
in the tree defining $X$ each vertex has at most one child which is not a
leaf. A $k$-lacunary metric space is a $k$-increasing metric space $X$ such
that in the tree defining $X$, each vertex has at most two children.
\end{definition}

Given integers $n, k$ and $\alpha>1$, we ask for the largest integer $m$ such
that every $n$-point metric space contains an $m$-point subspace which is
$\alpha$ embeddable in either an equilateral space or a $k$-lacunary space.
Formally, we define this quantity to be
$$D_k(\alpha,n)=R_{\left\{
\substack{ \text{$k$-lacunary spaces}\\ \text{or equilateral
spaces}}\right\}}(\alpha,n).$$ When $k>1$, this function exhibits a phase
transition at $\alpha=2$. When $k=1$, no phase transition occurs:

\begin{theorem}[The equilateral/lacunary dichotomy]\label{thm:4lacunary}
\mbox{\ }
\begin{enumerate}
\item For $\alpha>2$, $k>1$:
\[
c(\alpha,k) \cdot\frac{\log n}{\log \log n}\le D_k(\alpha,n) \le C(\alpha,k)
\cdot \frac{\log n}{\log \log n}.
\]
\item For $1<\alpha<2$, $k>1$:
\[
c(\alpha,k) \sqrt{\log n}  \le D_k(\alpha,n) \le C(k)\sqrt{\log n}.
\]
\item For any $\alpha>1$,
\[
c(\alpha) \log n\le D_1(\alpha,n)\le C\log n.
\]
\end{enumerate}
Here $c(\alpha,k),C(\alpha,k)>0$ depend only on $\alpha$ and $k$,
$c(\alpha)>0$ depends only on $\alpha$, $C(k)>0$ depends only on $k$, and
$C>0$ is an absolute constant.
\end{theorem}

\medskip

Dichotomic metric Ramsey problems have been studied for some time
now. The proof of Theorem~\ref{thm:bfm86} in \cite{bfm} uses
embedding into 1-increasing spaces (a class which contains both
$k$-lacunary spaces and equilateral spaces). Karloff, Rabani, and
Ravid \cite{krr} have studied the dichotomic problem in the
context of online computation, and obtained the lower bound above
for $D_{k}(4,n)$. In \cite{bbm} some of the upper bounds in
Theorem~\ref{thm:4lacunary} are proved.

\paragraph{Remark.} In graph theory, the study of
dichotomic Ramsey problems is usually not restricted to the
symmetric case. In our setting this translates to questions such
as: Given $\alpha\geq 1$, $k>1$, $e,f \in \mathbb{N}$, what is the
smallest $n$ such that every $n$-point metric space contains
either an $e$-point subspace which is $\alpha$ equivalent to an
equilateral space or an $f$-point subspace which is $\alpha$
equivalent to a $k$-lacunary space (respectively a binary
$k$-HST)? All our proofs extend in a straightforward manner to
give similarly tight bounds for the non-symmetric problems as
well.

\paragraph{Structure of the paper.} Our proof
of Theorem \ref{thm:second} relies on
Theorem~\ref{thm:ramsey-hst}. On the other hand, our proof for the
equilateral/lacunary dichotomy is elementary. We therefore start
with an elementary proof of Theorem~\ref{thm:4lacunary}
(Section~\ref{section:lacunary}), and then give a short proof,
based on a tool from \cite{blmn-phenomena}, of Theorem~\ref{thm:second}
(Section~\ref{sec:m-dichotomy}).

\section{The Equilateral/Lacunary Dichotomy}\label{section:lacunary}

In this section we prove Theorem~\ref{thm:4lacunary}.
A careful reading of the proof in \cite{bfm} shows that: Every $n$-point
metric space $M$ contains, for any $\epsilon>0$, a subspace $Y$ which is
$(1+\epsilon)$ embeddable in a $1$-increasing space, and
$|Y|=\Omega\left(\frac{\epsilon}{\log (1/\epsilon)}\log n \right)$.

We begin with a simplified proof of this, that works for $k$-increasing space
for any $k\geq 1$.

\begin{theorem} \label{thm:bfm} Fix an integer $n$, $0<\epsilon<1$ and $k\ge 1$.
Then any $n$ point metric space $X$ contains a subspace $Y\subset
X$ which is $(1+\epsilon)$ embeddable in a $k$-increasing space
and:
$$
|Y|\ge \frac{\epsilon}{6\log(12/\epsilon)\log(2k)}\cdot\log n.
$$
\end{theorem}

We start with the following simple lemma:

\begin{lemma}\label{lem:simple} Let $M$ be an $n$-point metric space and $0<\epsilon<1$.
Then there are $x  \in M$, $A  \subset M$ and $\lambda  \in [1,2]$ with the
following properties:

\medskip
\noindent{\bf 1)} $|A  |\ge \frac{\epsilon n}{4}$

\medskip
\noindent{\bf 2)} For every $z\in A  $, $\frac{\lambda
\diam(M)}{2(1+\epsilon)}\le d(x,z)<\frac{\lambda  \diam(M)}{2}$.

\end{lemma}

\begin{proof} Denote $\Delta=\diam(M)$. Take
$x,\bar{x}\in M$ such that $d(x,\bar{x})=\Delta$. The two sets
$Z=\{y\in M;\ d(y,x)<\Delta/2\}$ and $Z'=\{y\in M;\
d(y,\bar{x})<\Delta/2\}$ are disjoint, so we may assume that
$|Z|\ge n/2$. We split this set into layers.
$$
S_i=\left\{z\in M;\ (1 +\epsilon)^{i-1}\frac{\Delta}{2}\leq d(x,z)< (1 +
\epsilon)^i\frac{\Delta}{2} \right\}.
$$
Since $|Z|\ge n/2$, there exists $1\le i_0\le
\lceil\log_{1+\epsilon} 2\rceil$ such that
$$
|S_{i_0}|\ge \frac{n}{2 \lceil\log_{1+\epsilon} 2\rceil}\ge \frac{\epsilon
n}{4}.
$$

Take $A=S_{i_0}$ and $\lambda=\min\{2,(1+\epsilon)^{i_0}\}$ to obtain the
required result.
\end{proof}

We also need the following numerical lemma.
\begin{lemma} \label{lem:subseq}
Let $\{a_i\}_{i=1}^m$ a sequence of positive numbers,  satisfying
for any $i<j$, $a_j\leq 2 a_i$. Fix $\epsilon>0$ and $k\geq 1$.
Then there exists $L\subset\{1,\ldots,m\}$,
of cardinality $|L|\geq m/ (\lceil \log_{1+\epsilon}(2k) \rceil+1)$,
and a sequence $\{b_i\}_{i\in L}$ such
that for any
$i\in L$, $a_i\leq b_i \leq a_i (1+\epsilon)$, and for any
$i < j$ in $L$, either $b_{i}=b_j$ or $b_j\leq b_i/k$.
\end{lemma}
\begin{proof}
For every $a>0$, let $t(a)$ be the unique integer such that
$a\in\big((1+\epsilon)^{t(a)-1}, (1+\epsilon)^{t(a)}\big]$. Set
$r=\lceil \log_{1+\epsilon}(2k) \rceil+1$. For $j\in\{0,\ldots,
r-1\}$ define
\[ L_j=\{1\le i\le m; t(a_i) \equiv j\pmod r\}.\]
There is an integer $0\le j\le r-1$ such that $|L_j|\ge m/r$. Set
$L=L_j$. Define for $i\in L$, $b_i=(1+\epsilon)^{t(a_i)}$, hence
$a_i\leq b_i \leq a_i (1+\epsilon)$. Fix $i,j\in L$, $i<j$. Then
either $t(a_{i})=t(a_{j})$, in which case $b_i=b_j$, or otherwise
$t(a_{j})\not\in (t(a_{i})-r, t(a_{i})+r)$. We claim that
$t(a_{j})\leq t(a_{i_0})-r$. Indeed, otherwise $t(a_{j})\geq
t(a_{i})+r$, and therefore $a_{j}> a_{i} (1+\epsilon)^{r-1} \geq 2
a_{i}$, which  contradicts the assumptions. Therefore $b_{j}\leq
b_i/(1+\epsilon)^r\leq b_i/k$.
\end{proof}

\begin{proof}[Proof of Theorem~\ref{thm:bfm}] Set $\Delta=\diam(M)$. Let $x_1\in M$, $A_1\subset M$, $\lambda_1\in [1,2]$ be
as in Lemma~\ref{lem:simple}. Iterate this construction for $A_1$ until we
reach a singleton. We construct in this way $x_1,\ldots,x_m\in M$,
$\lambda_1,\ldots,\lambda_m\in [1,2]$ and $A_{m}\subset A_{m-1}\subset\ldots
A_1\subset M=A_0$ with the following properties:

\medskip
\noindent{\bf a)} $|A_{i+1}|\ge \frac{\epsilon}{4}|A_i|$.

\medskip
\noindent{\bf b)} For every $z\in A_{i+1}$, $\frac{\lambda_{i+1}
\diam(A_i)}{2(1+\epsilon)}\le d(x_{i+1},z)< \frac{\lambda_{i+1}
\diam(A_i)}{2}$

\medskip
\noindent{\bf c)} $A_m=\{x_m\}$ and $|A_{m-1}|>1$.

\medskip

These conditions imply in particular that $m\ge \frac{\log
n}{\log(4/\epsilon)}$.

The set $\{x_1,\ldots, x_m\}$ is therefore
$(1+\epsilon)$-equivalent to a metric similar to an increasing
space, but the labels are not monotonic. We solve this problem by
an appropriate sparsification. Put
$l_i=\frac{\lambda_{i}\diam(A_{i-1})}{2}$. Note that if $i>j$ then
$l_i\le 2l_j$. Indeed, this follows from the fact that $A_i\subset
A_j$ and $\lambda_i,\lambda_j\in [1,2]$. Apply
Lemma~\ref{lem:subseq} to the sequence $\{l_i\}_{i=1}^m$, and let
$\{b_i\}_{i \in L}$ be the resulting
sequence. Let $c_1> c_2>\cdots >c_s$ be such that $\{c_1,\ldots,
c_s\}= \{b_i;\  i\in L\}$.

For $i=1,\ldots, s$ define $J_i=\{h\in L ;\ b_h=c_i\}$ and put
$B_i=\cup_{h\in J_i}A_h$. Set also $B_0=M$. We construct a labelled tree $T$
as follows. The root of $T$ is $B_1$, and the rest of the vertices are
$\{x_i\}_{i\in L}$ and $\{B_i\}_{i=2}^s$. For $i\in L$, $x_i$ is a leaf of
$T$. The children of $B_i$  are $B_{i+1}$ and each of $\{x_h\}_{h\in J_i}$.
We label $T$ by setting for each $i\in L$, $\Delta(x_i)=0$, for $i=1,\ldots
s$ and $\Delta(B_i)=c_i$. By Lemma~\ref{lem:subseq}, $\Delta(B_{i+1})\le
\Delta(B_i)/k$.

Set $X=\{x_h\}_{h\in L_j}$. We have proved that $X$ is a $k$-increasing
space. Take $a,b\in L$, $a<b$. Assume that $a\in J_p$, $b\in J_q$, $p\le q$.
Since $x_b\in A_a$ we get that:
$$
d_X(x_a,x_b)=\Delta(B_p)=b_a\le
(1+\epsilon)l_a=\frac{(1+\epsilon)\lambda_{a}\diam(A_{a-1})}{2}\le
(1+\epsilon)^2 d(x_a,x_b),
$$
and
$$
d_X(x_a,x_b)=\Delta(B_p)=b_a\ge l_a=\frac{\lambda_a
\diam(A_{a-1})}{2} \ge d(x_a,x_b).
$$
This proves that $\{x_h\}_{h\in L}$ is $(1+\epsilon)^2$ equivalent to a
$k$-increasing space. The estimate on $|L|$ gives the required result.
\end{proof}

We can now deduce that any large metric space contains a large subspace which
is close to either an equilateral space or to a lacunary space.

\begin{proposition} \label{prop:krrHST2} Let $X$ be an $n$ point metric space, $k\ge 1$, and
$\epsilon>0$. Then $X$ contains a subspace $Y$ which is $(1+\epsilon)$
equivalent to either an equilateral space or a $k$-lacunary space, and such
that:
$$
|Y|\ge c\sqrt{\frac{\epsilon}{\log(2/\epsilon)}\cdot\frac{\log n}{\log k}},
$$
where $c$ is an absolute constant.
\end{proposition}

\begin{proof}
By Theorem~\ref{thm:bfm}, it is enough to prove that any $m$ point
$k$-increasing metric space contains (isometrically) either an equilateral
space of size $\sqrt{m}$ or a $k$-lacunary space of size $\sqrt{m}$.

Let $X$ be an $m$ point $k$-increasing space. Let $T$ be the tree defining
$X$. If $T$ contains an internal vertex with $\sqrt{m}$ leaves then these
leaves form an equilateral space. Otherwise $T$ contains at least $\sqrt{m}$
internal vertices which have at least one child which is a leaf. These leaves
form a $k$-lacunary metric space.
\end{proof}

For distortion $\alpha<2$, Proposition~\ref{prop:krrHST2} is tight. Here is a
matching upper bound.

\begin{proposition}\label{prop:ramseygraph}
For any $\alpha\in[1,2)$, any $k>1$ and any integer $n$ there exists an
$n$-point metric space $X$ such that no subset of $X$ with cardinality greater
that $\frac{c}{\log k}\sqrt{\log n}$ is $\alpha$ equivalent to an equilateral
space or a $k$-lacunary space. Here $c$ is an absolute constant.
\end{proposition}

The proof of Proposition~\ref{prop:ramseygraph} is based on the notion of
simple metric composition. This is a special
case of a more general definition that was introduced in \cite{blmn-phenomena}.
\begin{definition}[Simple Metric Composition] Let $M$, $N$
be two finite metric spaces and let $\beta\geq 1$. The $\beta$-composition of
$M$ and $N$ is a metric space on $M\times N$ which we denote by
$L=M_\beta[N]$. Distances in $L$ are defined by:
 \[ d_L((i,j),(k,l))= \begin{cases} d_N(j,l) & i=k \\
  \beta \gamma d_M(i,k) & i\neq k .\end{cases}  \]
where $\gamma=\frac{\diam(N)}{\min_{i\neq k} d_M(i,k)}$. It is easily checked
that the choice of the factor $\beta \cdot \gamma$ guarantees that $d_L$ is
indeed a metric.
\end{definition}

In words, first multiply the distances in $M$ by $\beta \cdot \gamma$, and
then replace each point of $M$ by an isometric copy of $N$.

We also use the notation $\Phi(X)=\frac{\diam(X)}{\min_{x,y\in
X,\, x\neq y} d_X(x,y)}$. This is the {\em aspect ratio}
of the metric space $X$, and in other words,
the Lipschitz distance between $X$ and
an equilateral space. We begin with three simple lemmas:

\begin{lemma}\label{lem:kbound} Let $X$ be a finite metric space which is
$\alpha$ embeddable in a $k$-lacunary metric space for some $k,\alpha>1$.
Then $ |X|\le 2+\log_{k}\left(\alpha \Phi(X)\right). $
\end{lemma}
\begin{proof}
Let $Y$ be a $k$-lacunary space that is $\alpha$ equivalent to $X$. Hence
$\Phi(Y)\leq \alpha \Phi(X)$. A simple induction on $|Y|$ shows that for any
$k$-lacunary space $Y$, $\Phi(Y)\geq k^{|Y|-2}$.
\end{proof}

\begin{lemma}\label{lem:unifcomb} Let $M$, $N$ be finite metric spaces,
and let $\beta> \alpha\ge1$.
Then every subset $S\subset M_\beta[N]$ with $\Phi(S)\leq \alpha$
is 1-embeddable either in $M$ or in $N$.
\end{lemma}

\begin{proof}
For every $x\in M$ denote $D_x=\{(x,y)\in M\times N;\ y\in N\}$. If $S\subset
D_x$ for some $x\in M$ then $S$ is $1$-embeddable in $N$. If for each $x\in
M$, $|S\cap D_x|\le 1$ then $S$ is $1$-embeddable in $M$. Otherwise there are
$a,b,c\in S$ and $x, y\in M$, $a\neq b$, $x\neq y$, such that $a,b\in D_x$
and $c\in D_y$. Hence:
$$
\frac{d_S(a,c)}{d_S(a,b)}\ge \frac{\beta\gamma\min_{u\neq
v}d_M(u,v)}{\max_{u,v\in N}d_N(u,v)}\ge {\beta}>\alpha,
$$
which contradicts the fact that $\Phi(S)\leq \alpha$.
\end{proof}

\begin{lemma}\label{lem:lacunarycomb} Let $M$, $N$ be finite metric spaces. Fix
$k,\alpha\ge 1$ and $\beta\ge \max\left\{1,\frac{\alpha}{k}\right\}$. Then
every $S\subset M_\beta[N]$ which is $\alpha$-embeddable in a $k$-lacunary
metric space has a subset $T\subset S$ that is $1$-embeddable in $M$ and
$S\setminus T$ is $1$-embeddable in $N$.

\end{lemma}

\begin{proof} Let $a_1\ge\ldots\ge a_n>0$ be a $k$-lacunary sequence, i.e.,
$a_{i+1}\le a_i/k$. It is easy to verify that for every four distinct
integers $1\le i_1,i_2,i_3,i_4,\le n$,
$$
\max\{a_{\min\{i_1,i_2\}},a_{\min\{i_3,i_4\}}\}\ge
k\min\{a_{\min\{i_1,i_3\}}, a_{\min\{i_1,i_4\}},
a_{\min\{i_2,i_3\}},a_{\min\{i_2,i_4\}}\}.
$$
Since $S$ is $\alpha$-embeddable in a $k$-lacunary space,
it follows that for every distinct $x_1,x_2,x_3,x_4\in S$,
\begin{eqnarray}\label{eq:fourpoint}
\max\{d_S(x_1,x_2),d_S(x_3,x_4)\}\ge \frac{k}{\alpha}
\min\{d_S(x_1,x_3),d_S(x_1,x_4),d_S(x_2,x_3),d_S(x_2,x_4)\}.
\end{eqnarray}

As before, denote $\gamma=\frac{\max_{x,y\in N}d_{N}(x,y)}{\min_{x\neq y}
d_{M}(x,y)}$ and for $x\in M$, $D_x=\{(x,y)\in M\times N;\ y\in N\}$. It is
sufficient to prove that there is at most one $x\in M$ such that $|S\cap
D_x|>1$. This is true since otherwise there would be four distinct points
$p,q,r,s\in N$ and two distinct points $x,y\in M$ such that $p,q\in D_x$ and
$r,s\in D_y$. Now:
\[
\frac{\max\{d_S(p,q),d_S(r,s)\}}{\min\{d(p,r),d(p,s),d(q,r),d(q,s)\}}\le
\frac{\max_{u,v\in N}d_N(u,v)}{\beta\gamma\min _{u\neq
v}d_M(u,v)}=\frac{1}{\beta}< \frac{k}{\alpha},
\]
which contradicts (\ref{eq:fourpoint}).
\end{proof}

\begin{proof}[Proof of Proposition~\ref{prop:ramseygraph}]
Let $G=(V,E)$ be a graph of diameter 2 on
$\left\lceil2^{\sqrt{\log n}}\right\rceil$ vertices, with
no independent sets and no cliques larger than
$C \sqrt{\log n}$, here $C$ is an absolute constant. It is
well-known and easy to prove that almost all
graphs have these properties, see \cite{erdos, bollobasext}.

Let $M$ be the metric defined by $G$.
Define $M_1=M$, and $M_i=M_\beta[ M_{i-1}]$, where $\beta=2$. First we
prove by induction that for each $i\ge 1$, if $S\subset M_i$ is $\alpha$
embeddable in an equilateral space then $|S|\le C\sqrt{\log n}$. For
$i=1$ consider a subset $S\subset M$ that is $\alpha$-embeddable in an
equilateral space. Since $\alpha<2$, and $G$ has diameter 2, all
the distances in $S$ must be either $1$ or $2$. Thus $S$ is either a clique
or an independent set, so that $|S|\le C\sqrt{\log n}$. Now let $i>1$ and
consider $S\subset M_i=M_\beta[M_{i-1}]$
that is $\alpha$-embeddable in an equilateral
space. By Lemma~\ref{lem:unifcomb}, $S$ is $1$-embeddable in either $M$ or
$M_{i-1}$, which by induction implies that $|S|\leq C \sqrt{\log n}$.

We now prove by induction on $i$ that if $S\subset M_i$ is
$\alpha$-embeddable in a $k$-lacunary space then $|S|\le
i\left(1+\log_k(2\alpha)\right)$. For $i=1$ this follows from
Lemma~\ref{lem:kbound}. For $i>1$ let $S\subset M_i=M_\beta[M_{i-1}]$
be $\alpha$-embeddable in a $k$-lacunary space. By
Lemma~\ref{lem:lacunarycomb} this implies that there is $A\subset S$
that $1$-embeds into $M$ and $S\setminus A$ $1$-embeds into $M_{i-1}$. By
the induction hypothesis:
$$
|S|=|A|+|S\setminus A|\le
1+\log_k(2\alpha)+(i-1)\left(1+\log_k(2\alpha)\right)=
i\left(1+\log_k(2\alpha)\right).
$$
Take $t=\lceil \sqrt{\log n}\rceil$ and note that $|M_t|\geq n$.
The space $X=M_t$ satisfies our claim.
\end{proof}

For every $b > a > 1$,
it is easy to extract from every $a$-lacunary sequence
a long $b$-lacunary subsequence, by skipping each time
appropriately many terms in the sequence.
We record this simple fact for future reference.

\begin{lemma}\label{lem:subsequence} For every $b>a>1$,
every $a$-lacunary sequence of
length $n$ contains a subsequence of length $\frac{n}{\lceil 1+\log_a
b\rceil}$ which is $b$-lacunary. Hence, any $n$ point $a$-lacunary metric
space contains a $b$-lacunary subspace of the above size.
\end{lemma}

Using a technique similar to \cite{krr}, we now prove:

\begin{proposition}\label{prop:lowerklac} For any $k>1$, $\alpha>2$
and any integer $n$, every $n$ point metric space contains a subspace of
cardinality at least $\frac{\log(\alpha/2)}{2\log (\alpha k)}\cdot\frac{\log
n}{\log \log n}$ which is $\alpha$-embeddable in either an equilateral space
or a $k$-lacunary space.
\end{proposition}
\begin{proof}

Let $(M,d)$ be an $n$-point metric space. Denote
$\Delta=\diam(M)$, and let $x,\bar{x}\in M$ be a diametrical pair,
i.e., $d(x,\bar{x})=\Delta$. Let $x=x_1,\ldots, x_s$ be a maximal
subset in $M$ containing $x$ such that for every $i\neq j$,
$d(x_i,x_j)\geq \Delta/\alpha$. Clearly $\{x_1,\ldots x_s\}$ is
$\alpha$-equivalent to an equilateral space, so that if $s\ge \log
n$ we are done. As usual, we denote by $B(x,r) = \{y \in M ;\
d(x,y) < r\}$ the open ball of radius $r$ around $x$. Let
$$
A_1=B\left(x_1,\frac{\Delta}{\alpha}\right) \quad \mathrm{and}
\quad A_{i+1}=
B\left(x_{i+1},\frac{\Delta}{\alpha}\right)\setminus\bigcup_{j=1}^iB\left(x_j,\frac{\Delta}{\alpha}\right).
$$
Assume that $s<\log n$. Since $\cup_{i=1}^sA_i=
\cup_{i=1}^sB(x_i,\Delta/\alpha)=M$, it follows that there exists
$1\le i\le n$ such that $|A_i|\ge n/\log n$. Observe that there
exists $y\in M$ such that $d(y,A_i)\ge \Delta/\alpha$. Indeed, if
$i>1$ then $A_i\subset M\setminus B(x,\Delta/\alpha)$ so that we
can take $y=x$. Otherwise, since $\alpha>2$, for every $z\in
A_1=B(x,\Delta/\alpha)$,
$$
d(z,\bar{x})\ge d(\bar{x},x)-d(x,z)\ge
\Delta-\frac{\Delta}{\alpha}>\frac{\Delta}{\alpha},
$$
so that we can take $y=\bar{x}$. Note also that
$
\diam(A_i)\le \diam(B(x_i,\Delta/\alpha))\le
\frac{2}{\alpha}\Delta$.

Iterating this construction we get a sequence of points
$z_1,\ldots, z_m\in M$, and a decreasing sequence of subsets
$\{z_m\}=F_m\subset F_{m-1}\subset\ldots \subset F_1\subset F_0=M$
such that for each $i$, $z_i\in F_{i-1}$, $d(z_i,F_i) \ge
\diam(F_{i-1})/\alpha$, $\diam(F_i)\le
\frac{2}{\alpha}\diam(F_{i-1})$ and $|F_i|\ge |F_{i-1}|/\log n$.
By induction, $|F_i|\ge n/(\log n)^i$ and since $|F_m|=1$,
necessarily $m\ge \log n/\log\log n$. Moreover, the sequence
$\{\diam(F_i)\}_{i=0}^{m-1}$ is $(\alpha/2)$-lacunary and for
$1\le i<j\le m$,
$$
\frac{\diam(F_{i-1})}{\alpha}\le d(z_i,F_i)\le d(z_i,z_j)\le \diam(F_{i-1}).
$$
This proves that $\{z_i,\ldots,z_m\}$ is $\alpha$ equivalent to a
$(\alpha/2)$-lacunary metric space. If $k<\alpha/2$ we are done. Otherwise,
we can apply Lemma~\ref{lem:subsequence} to find a subset of
$\{z_1,\ldots,z_m\}$ which is $\alpha$ embeddable in a $k$-lacunary space and
with cardinality at least:
$$
\frac{m}{\lceil 1  +\log_{\alpha/2}k\rceil}\ge \frac{\log(\alpha/2)}{2\log
(\alpha k)}\cdot \frac{\log n}{\log \log n}.
$$
\end{proof}

We can now establish the equilateral/lacunary dichotomy:

\begin{proof}[Proof of Theorem~\ref{thm:4lacunary}.]
The lower bound in part ${\bf 1)}$ was proved in
Proposition~\ref{prop:lowerklac}.
In \cite{bbm}, Proposition 29 it is proved that for
$2<\alpha<k$,
$$
D_{\alpha}(n,k)\le C \frac{\log \alpha}{\log k}\cdot \frac{\log n}{\log \log
n}.
$$
In the case $\alpha\ge k$, let $M$ be a metric space and $N\subset M$ a subset
which is $\alpha$ embeddable in either an equilateral space or a $k$-lacunary
space. Apply Lemma~\ref{lem:subsequence} with $b=\alpha^2$ and $a=k$. We
deduce that there is $N'\subset N$ which is $\alpha$ embeddable in either an
equilateral space or a $\alpha^2$-lacunary space such that $|N'|\ge
\frac{|N|\log k}{2\log \alpha}$. By the above stated result from \cite{bbm},
$$
\frac{|N|\log k}{2\log \alpha}\le |N'|\le C \frac{\log \alpha}{\log
(\alpha^2)}\cdot \frac{\log n}{\log \log n},
$$
which implies the required result.

Part ${\bf b)}$ is a combination of Proposition~\ref{prop:krrHST2} and
Proposition~\ref{prop:ramseygraph}. Part ${\bf c)}$ is a combination of
Theorem~\ref{thm:bfm} and part 4 of Proposition 29 in \cite{bbm}.
\end{proof}

\section{The Metric Dichotomy} \label{sec:m-dichotomy}

Our main aim in this section is to prove Theorem~\ref{thm:second}. We begin,
however, with a proof of Proposition~\ref{prop:approx-equiv}.

\begin{proof}[Proof of Proposition~\ref{prop:approx-equiv}]
To prove the first assertion note that if
$x,y,z$ are three distinct leaves in a binary tree $T$,
then there are $p,q\in T$ such that $p\neq q$,
$q$ is a descendant of $p$ and $\{\lca(x,y), \lca(x,z), \lca(y,z)\}
=\{p,q\}$. Since $\Delta(p)/\Delta(q)\ge k$, the Lipschitz distance between
the triangle $\{x,y,z\}$ and an equilateral triangle is at least $k$.

To prove the second assertion, denote $\Delta=\diam(M)$, and let
$x,\bar{x}\in M$ be a diametrical pair: $d(x,\bar{x})=\Delta$. Let
$B_x=B(x,\Delta/k )$, and $B_{\bar{x}}=B(\bar{x},\Delta/k)$. Since
$k>2$, the triangle inequality implies that $B_x \cap
B_{\bar{x}}=\emptyset$. We claim that $B_x \cup B_{\bar{x}}=M$.
Otherwise, there is a $ y\in M$ such that $d(x,y) \ge \Delta/k$,
and $d(\bar{x},y) \ge \Delta/k$. But this implies that
$x,\bar{x},y$ are three points for which $c_{\{\text{equilateral\
spaces}\}}(\{x,\bar{x},y\}) \le k$, contrary to our assumption.

We proceed to construct a binary $\frac{k}{2}$-HST $L$ with
$\diam(L)=\diam(M)$, and a non-contractive embedding
$g:M\hookrightarrow L$ with $\|g\|_{\text{lip}}\le \frac{k}{k-2}$.
Inductively, assume we have already constructed $g_1:B_{{x}}
\hookrightarrow L_1$, $g_2:B_{\bar{x}} \hookrightarrow L_2$, where
$L_1$, $L_2$ are binary ${k}$-HSTs with $\diam(L_1)=
\diam(B_{x})$, $\diam(L_2)= \diam(B_{\bar{x}})$, and $g_1,g_2$ are
non-contractive embeddings with $\|g_1\|_{\text{lip}}\leq
\frac{k}{k-2}$, and $\|g_2\|_{\text{lip}}\leq \frac{k}{k-2}$.
Define $L=L_1\cup L_2$, and $g:M \hookrightarrow L$ by
$g|_{B_{{x}}}=g_1$, $g|_{B_{\bar{x}}}=g_2$. Set the distance
between any point in $L_1$ and any point in $L_2$ to be $\Delta$.
Since $\max\{\diam(L_1), \diam(L_2)\}\leq  \frac{2\Delta}{k}$, $L$
is a binary $\frac{k}{2}$-HST, $g$ is non-contractive,
$\diam(L)=\diam(M)$, and
 \[ \|g\|_{\text{lip}}\leq
\max \Bigl\{ \|g_1\|_{\text{lip}}, \|g_2\|_{\text{lip}},
\frac{\Delta}{\Delta- 2\Delta/k} \Bigr \} \leq \frac{k}{k-2}.
 \]
\end{proof}

We begin with the upper bounds on $E_k$ and $F_k$ for distortions smaller
than 2. We give both a bound for $E_k$ and for $F_k$, since the bound for
$E_k$ holds for any $k> 1$, whereas for $F_k$, it holds only for $k>2$.

\begin{lemma}\label{lem:kbound2} Let $X$ be a finite metric space which is
$\alpha$-embeddable in a binary $k$-HST for some $k,\alpha>1$. Then $
|X|\le 2^{1+\log_{k}\left(\alpha \Phi(X)\right)}. $
\end{lemma}
\begin{proof}
Let $Y$ be a binary $k$-HST that is $\alpha$ equivalent to $X$. Hence
$\Phi(Y)\leq \alpha \Phi(X)$. The tree defining $Y$ is binary and
its depth is therefore
$\ge \log_2 |Y|$. A simple induction on $|Y|$ proves
that for any binary $k$-HST $Y$, $\Phi(Y)\geq k^{\log_2|Y|-1}$.
\end{proof}

\begin{proposition}\label{prop:sec-ub-1}
For any $\alpha\in[1,2)$, any $k>1$ and any integer $n$ there exists an $n$
point metric space $X$ such that no subset of $X$ with cardinality greater
than $\frac{c}{\log k} \frac{\log n}{\log \log n}$ is $\alpha$-equivalent to
an equilateral space or a binary $k$-HST. Here $c$ is an absolute constant.
\end{proposition}
\begin{proof}
Again we recall that almost every graph on
$s=\left\lceil2^{2(1+\log_k(2\alpha)) \frac{\log n}{\log\log n} }\right\rceil$
vertices has diameter 2 and its
independence number and clique number are $\le C \log s $, for some
absolute constant $C$.

Let $G$ be such a graph and let $M$ be its metric.
Next, define $M_0=\{a\}$, and $M_i=M_\beta[ M_{i-1}]$, where $\beta=2$.

Similarly to the proof of
Proposition~\ref{prop:ramseygraph}, for each $i\ge 1$, if $S\subset M_i$ is
$\alpha$-embeddable in an equilateral space then $|S|\le C\log s $.
For $i=1$, if $S\subset M$ is $\alpha$-embeddable in an
equilateral space, then $S$ must either be a clique
or an independent set, since $\alpha<2$.
Consequently, $|S|\le C\log s $. Now let
$S\subset M_i=M_\beta[M_{i-1}]$ be $\alpha$-embeddable in an equilateral
space for some $i >1$. By Lemma~\ref{lem:unifcomb},
$S$ is $1$-embeddable in either $M$ or
$M_{i-1}$, which by induction implies that $|S|\leq C \log s$.

We now prove by induction on $i$ that if $S\subset M_i$ is
$\alpha$-embeddable in a binary $k$-HST then $|S|\le
2^{i(1+\log_k(2\alpha))}$. For $i=0$ this is obvious. Assume that
$i>0$ and $S\subset M_i=M_\beta[M_{i-1}]$ is $\alpha$ embeddable
in a binary $k$-HST. Partition $S$ to $S=S_1\cup\ldots \cup
S_\ell$ such that each $S_j$ is a subset of a different ``copy" of
$M_{i-1}$. Note that $S_j\subset M_{i-1}$  is  $\alpha$ embeddable
in a binary $k$-HST, and by the inductive hypothesis $|S_j| \leq
2^{(i-1)(1+\log_k(2\alpha))}$. Pick a representative from each
$S_j$, and denote the set of representatives by $S'$, $|S'|=\ell$.
As $S'\subset S$ it is also $\alpha$ embeddable in a binary
$k$-HST (the defining tree is a subtree of the tree defining the
HST of $S$). The metric of $S'$ is a dilation of a subset of $M$,
and so by Lemma~\ref{lem:kbound2}, $\ell\leq
2^{1+\log_k(2\alpha)}$. We can therefore estimate
 \[ |S|\leq \ell \max_{1\le j\le \ell}|S_j|  \leq 2^{1+\log_k(2\alpha)}
2^{(i-1)(1+\log_k(2\alpha))}=2^{i(1+\log_k(2\alpha))}.\]

Note that $|M_t|\geq n$
for $t=\lceil \frac{\log \log n}{2(1+\log_k(2\alpha))}\rceil$.
The space $X=M_t$ satisfies the Proposition.
\end{proof}

We give a similar upper bound for the equilateral/triangular
variant of this dichotomy.

\begin{proposition}\label{prop:sec-ub-2}
For any $\alpha\in[1,2)$, any $k>2$ and any integer $n$ there
exists an $n$-point metric space $M$ such that no subset of $M$
with cardinality greater than $c \log_{k/2} \frac{k}{k-2} \cdot
\frac{\log n}{\log \log n}$ is $\alpha$-equivalent to an
equilateral space or a space in which no triangle is
$k$-equivalent to an equilateral triangle. Here $c$ is an absolute
constant.
\end{proposition}
\begin{proof}
The proof is almost identical to the proof of
Proposition~\ref{prop:sec-ub-1}. The only change is the reference to
Lemma~\ref{lem:kbound2}. Here instead we use the following claim:
Any finite metric
space $X$ which is $\alpha$ embeddable in space in which no triangle is
$\le k$-equilateral, for some $k>2$, $ \alpha>1$
satisfies $|X|\le 2^{1+\log_{k/2}\left(\alpha \frac{k}{k-2} \Phi(X)\right)}.
$ This fact is an immediate consequence of Lemma~\ref{lem:kbound2} and
Proposition~\ref{prop:approx-equiv}.
\end{proof}

The proofs of the lower bounds on $E_k$ use the following simple structural
lemma.
\begin{lemma} \label{lem:large-b-HST}
Let $T$ be a rooted tree with $n$ leaves, in which each vertex has
at most $h\ge 2$ children. Then $T$ contains a binary subtree with
at least $n^{1/\log_2h}$ leaves.
\end{lemma}
\begin{proof}
By induction on the size of $T$. Let $h'\leq h$
be the number of children of $T$'s root $r$.
Let $T_i$ be the subtree rooted at $r$-th $i$-th child
and let $n_i$ be the number of leaves in $T_i$, where
$n_1\geq n_2 \geq \ldots \geq n_{h'}$ and $\sum_{i=1}^{h'} n_i=n$.
By the induction hypothesis, $T_i$ has a binary subtree with at
least $n_i^{1/ \log_2 h}$ leaves. We form a binary subtree of $T$
by joining the binary subtrees of $T_1$ and $T_2$. Together they
have at least $n_1^{1/\log_2 h} + n_2^{1/\log_2 h}$ leaves,
which is $\ge n^{1/\log_2 h}$ as we now show. First,
$n_1^{1/\log_2 h} + n_2^{1/\log_2 h}
\ge 2 (\frac{n_1+n_2}{2})^{1/\log_2 h}$ since the function
$f(x)=x^{1/\log_2 h}$ is concave. Also, $f$ is increasing,
and $\frac{n_1+n_2}{2} \ge \frac{n}{h'} \ge \frac {n}{h}$.
Consequently, $n_1^{1/\log_2 h} + n_2^{1/\log_2 h}
\ge 2 (\frac{n}{h})^{1/\log_2 h} = n^{1/\log_2 h}$, as claimed.
\end{proof}

The following is a short argument proving the lower bound on $E_k$ for
distortions larger than 2.

\begin{proposition}\label{prop:sec-lb-2}
For any $\e\in(0,1)$, and $k\geq 1$,
\[ E_k(2+\e,n) \geq \exp \left ( \sqrt{\frac{c \e}{\log (2k/\e)} \log n} \right ). \]
\end{proposition}
\begin{proof}
Let $M$ be an arbitrary $n$-point metric space. By
Theorem~\ref{thm:ramsey-hst}, it contains a subset $N\subset M$ that is
$(2+\e)$-equivalent to a $k$-HST and $|N| \geq n^{\frac{c\e}{\log (2k/\e)}}$.
Let $T$ be the tree defining this $k$-HST. The claim is now proved by taking
either a large equilateral subspace of $T$ or a large binary subtree of $T$
according to Lemma~\ref{lem:large-b-HST}, where $h=\exp \bigl ( \sqrt{\frac{c
\e}{\log (2k/\e)} \log n} \bigr )$.

\end{proof}

We note that the above proposition can also be proved by arguments
similar to those from \cite{bkrs}.

\medskip

The lower bound on $E_k$ for distortions smaller than 2 is only slightly more
complicated:

\begin{proposition} \label{prop:sec-lb-1}
For any $>\e>0$, $k\geq 1$,
\[ E_k(1+\e,n)\geq \frac{c \e}{\log (1/\e) \log (k/\e)}\cdot \frac{\log n}{\log \log n},
\]
where $c$ is a universal constant.
\end{proposition}
\begin{proof}
Essentially, we repeat the argument from
Proposition~\ref{prop:sec-lb-2}, and find either an $h$-point
subspace that is 3-equivalent to equilateral space or an
$n^{1/\log h}$-point subspace that is 3 equivalent to binary
$k$-HST. In order to improve the distortion to $1+\e$,  we invoke
in the first case the classical Ramsey theorem to find a $\ge \log
h$-point subspace which is $(1+\e)$-equivalent to an equilateral
space, whereas in the second case we observe that by optimizing
the distances in the binary $k$-HST, we improve the distortion. We
choose $h\approx n^{1/\log\log n}$, so that $\log h \approx
n^{1/\log h}$. We now turn to the actual arguments.

Let $M$ be an $n$-point metric space. Denote $k'=\max\{k,2+2/\e\}$. By
Theorem~\ref{thm:ramsey-hst}, $M$ contains a subspace $N\subset M$ that is
$3$-equivalent to a $(3k')$-HST, $X$, via a non-contractive bijection $f:N\to
X$, and $|N| \geq n^{\frac{c}{\log k'}}=s$. Let $h= s^{1/ \log \log n}$.
Denote by $T$ the tree defining $X$. We distinguish between two cases.

\paragraph{Case 1.} $T$ has a vertex $u$ with out-degree exceeding $h$.
Let $v_0,\ldots, v_h$ be distinct children of $u$. For each $0\le i\le h$
take $x_i\in N$ such that $f(x_i)$ is a leaf of $T$ which is a descendant of
$v_i$. For every $0\le i<j\le h$, $d(x_i,x_j)\in [\Delta(u)/3, \Delta(u)]$,
so that there is a unique integer $c(i,j)\in \{1,2,\ldots,\lfloor
\log_{1+\e}3\rfloor\}$ for which:
$$
d(x_i,x_j)\in
\left[\frac{\Delta(u)}{(1+\e)^{c(i,j)+1}},\frac{\Delta(u)}{(1+\e)^{c(i,j)}}\right).
$$
Set $D=\lfloor \log_{1+\e}3\rfloor$ and color the edges of the complete graph
on $\{0,\ldots,h\}$ by assigning the color $c(i,j)$ to the edge $[i,j]$. By
the classical Ramsey theorem there is a subset $N'\subset \{x_1,\ldots,x_h\}$
of size at least $\frac{\log h}{D \log D}\ge \frac{c \e}{\log (1/\e)}
\frac{\log s}{\log\log n}$ on which the induced complete subgraph is
monochromatic.
This subset is $(1+\e)$-equivalent to an equilateral space.

\paragraph{Case 2.} All the vertices in $T$ have out-degree at most $h$. In this case, by
Lemma~\ref{lem:large-b-HST}, $T$ contains a binary subtree $S$
with at least $s^{1/ \log_2 h}= \log n$ leaves. Set $L=f^{-1}(S)$.
Then $|L|=|S| \geq \log n$ and $L$ is 3-equivalent to a binary
$(3k')$-HST $S$. In order to improve the distortion we
\emph{change the labels of $S$}. Denote by $\Delta(\cdot)$ the
original labels on $S$ (inherited from $T$). We define new labels
$\Delta'(\cdot)$ on $S$ as follows. For each vertex $u\in S$,
denote by $T_1$ and $T_2$ the subtrees rooted at $u$'s children.
We define $\Delta'(u)=\max\{d_M(x,y);\  x\in f^{-1}(T_1), \, y\in
f^{-1}(T_2)\}$, and claim that the resulting labelled tree is a
binary $k'$-HST which is $\frac{k'}{k'-2}$ equivalent to $L$.
Indeed, let $u,v\in S$ with $v$ a child of $u$. Since the
distances in $(S,\Delta)$ are larger than the distances in $M$ by
a factor at most 3, $\Delta(u)/3\leq \Delta'(u)$. On the other
hand, since $\Delta$ defines $3k'$-HST, $\Delta'(v)\leq
\Delta(v)\leq \Delta(u)/(3k')$, so that the resulting tree
$(S,\Delta')$ is indeed $k$-HST. To bound the distortion, let
$x,y$ be two distinct points in $L$. So $f(x)$, $f(y)$ are
distinct leaves of $S$ and assume that $\lca(f(x),f(y))=u$,
$f(x)\in T_1$, $f(y)\in T_2$, where $T_1$ and $T_2$ are subtrees
rooted at children of $u$. Then $d_M(x,y)\leq \Delta'(u)$. On the
other hand fix $a\in f^{-1}(T_1)$ and $b\in f^{-1}(T_2)$ for which
$d_M(a,b)=\Delta'(u)$. Then
\begin{eqnarray*}
d_L(x,y)&\geq&
d_L(a,b)-d_L(a,x)-d_L(b,y)\\
&\ge& \Delta'(u)- \Delta(\lca(f(a),f(x)))-\Delta(\lca(f(b),f(y)))\\
&\ge& \Delta'(u)-2\frac{\Delta(u)}{3k'}\\
&\ge& \Delta'(u) - \frac{2\Delta'(u)}{k'}=\frac{k'-2}{k'}\Delta'(u).
\end{eqnarray*}
since $k'\geq k$, and $\frac{k'}{k'-2} \leq 1+\e$, $L$ is $(1+\e)$-equivalent
to binary $k$-HST $(S,\Delta')$.
\end{proof}

\begin{proof}[Proof of Theorem~\ref{thm:second}]
The lower bound for $E_k(\alpha,n)$, $\alpha>2$, is contained in
Proposition~\ref{prop:sec-lb-2}. The upper bound for $F_k(\alpha,n)$,
$\alpha,k>2$ can be derived from results of \cite{bbm}, where it is proved
that for $1<\alpha<k$, $E_k(\alpha,n)\leq 2^{2 \sqrt{\frac{\log \alpha}{\log
k} \log n}}$. In order to prove the upper bound on $E_k(\alpha,n)$ for
$1<k\leq \alpha$, we use another lemma from~\cite{bbm}: for any $k>1$ and any
$h\in\mathbb{N}$, any $n$-point $k$-HST contains isometrically a subspace of
size  $n^{1/h}$ which is a $k^h$-HST. It is easy to observe that if we
start with a \emph{binary} $k$-HST the resulting subspace is a
\emph{binary} $k^h$-HST. Therefore, if $1<k\leq \alpha$, we take $h= \lfloor
1+\log _k \alpha\rfloor$. Using the discussion above, and the fact $h\geq 1$,
we deduce that
\[ [E_k(\alpha,n)]^{1/h} \leq E_{k^h}(\alpha,n)\leq 2^{2 \sqrt{\frac{\log \alpha}{h\log k} \log n}}. \]
As $\frac{\log \alpha}{h \log k} \leq 1$, we conclude that $E_k(\alpha,n)\leq
2^{2h\sqrt{\log n}}$. For $\alpha, k>2$, by
Proposition~\ref{prop:approx-equiv},
 \[ F_k(\alpha,n)\leq E_{k/2}(\alpha \tfrac{k}{k-2},n)\leq
 2^{2\left[1+\log_{k/2}\left(\tfrac{k}{k-2} \alpha\right)\right]\sqrt{\log n}}. \]

The lower bound for $E_k(\alpha,n)$, $\alpha\in(1,2)$, is contained in
Proposition~\ref{prop:sec-lb-1}. The upper bound for $F_k(\alpha,n)$ $k>2$,
$\alpha\in(1,2)$ is contained in Proposition~\ref{prop:sec-ub-2}. The
extension of this upper bound for $E_k(\alpha,n)$, $k\geq 1$,
$\alpha\in(1,2)$, is contained in Proposition~\ref{prop:sec-ub-1}.
\end{proof}

\bibliographystyle{plain}
\bibliography{dichotomy}

\bigskip
\bigskip

\noindent Yair Bartal, Institute of Computer Science, Hebrew
University, Jerusalem 91904, Israel. \\{\bf yair@cs.huji.ac.il}

\medskip
\noindent Nathan Linial, Institute of Computer Science, Hebrew
University, Jerusalem 91904, Israel. \\ {\bf nati@cs.huji.ac.il}

\medskip
\noindent Manor Mendel, Institute of Computer Science, Hebrew
University, Jerusalem 91904, Israel. \\ {\bf
mendelma@cs.huji.ac.il}

\medskip
\noindent Assaf Naor, Theory Group, Microsoft Research, One
Microsoft Way 113/2131, Redmond WA 98052-6399, USA. \\ {\bf
anaor@microsoft.com}

\bigskip
\bigskip
\noindent 2000 AMS Mathematics Subject Classification: 52C45,
05C55, 54E40, 05C12, 54E40.
\end{document}